\documentclass[12pt]{article}
\usepackage{amssymb}
\def \R{I \! \! R}
\def \N{I \! \! N}

\newcommand{\1}{1 \!\! 1}
\newtheorem{Theorem}{Theorem}
\newtheorem{Corollary}[Theorem]{Corollary}

\newtheorem{Proposition}[Theorem]{Proposition}

\newcommand{\F}{{\cal F}}

\newcommand{\A}{{\cal A}}

\begin{document}
\title{The Quicksort process}
\author{Mahmoud Ragab\\Mathematics Department\\
Faculty of Science\\
Al Azhar University\\
Nasr City,11884\\
Cairo\\Egypt
\\and \\Uwe Roesler\\Mathematisches Seminar\\
Christian-Albrechts-Universit\"at zu Kiel\\Ludewig-Meyn-Strasse 4
\\24098 Kiel\\Germany} \maketitle

\begin{abstract}
Quicksort on the fly returns the input of $n$ reals in increasing natural 
order during the sorting process. Correctly normalized the running time up 
to returning the $l$-th smallest out of $n$ seen as a process in $l$ converges 
weakly to a limiting process with path in the space of cadlag functions. 
\end{abstract}

AMS-classification: Primary 60 F 05, Secondary 68 P 10, 60 K 99

{\it Keywords and phrases:\/} sorting, Quicksort, divide and conquer
algorithm, random algorithm, running time analysis, stochastic process,
cadlag functions, asymptotics, Skorodhod metric

\section{Introduction}
Quicksort was chosen as one of the 10 most important algorithms. Quicksort serves also as 
a challenging random divide-and-conquer algorithm for a mathematical analysis. Starting 
with the worst case, the best case and the expected running time \cite{Knuth} we know 
nowadays much finer results on the limiting distribution, the existence \cite{Regnier} 
via martingale methods and a characterization \cite{Roesler91} as a stochastic fixed point.   
The running time for many versions, actually all versions I know of, can be analyzed \cite{MarRoe,MarRou} 
by the tools (contraction method and Weighted Branching Process \cite{RoeslerRueschendorf}) invented 
for Quicksort. 

Around 2004 Conrado Mart\'{i}nez \cite{Mar} came up with the algorithm Partial Quicksort, a mixture 
of Quickselect and standard Quicksort. Given as input a list of $n$ numbers, find the numbers 
between the $k$-th and the $l$-th smallest in order. A special case is Quicksort on the fly, 
where the algorithm tries to finish always the left most list. The output is 
first the smallest then the second smallest element and so on as an incoming stream. 

What can we say about the running time up to see the $l$-th smallest? Mart\'{i}nez gave an explicit 
formula for the expected time (\ref{1.2a}) and showed an asymptotic result. This was an indication 
for a possible distributional result. In 2010 Mart\'{i}nez and Roesler showed convergence of one dimensional 
marginals to a limit, which required already some effort. This paper is on the convergence of finite 
dimensional marginals and, moreover, the existence of a limiting process for Quicksort on the fly 
with cadlag paths in $D.$ Most of this work is in the dissertation of Ragab \cite{Ragab}.  

The existence of the Quicksort process, Theorem \ref{MainTheorem}, uses a specific construction via 
a Weighted Branching Process (WBP), see section 2 for the WBP and section 3 for the existence. The 
construction uses a specific sequence of rvs ($R_n$) with values in $D$ and the convergence is a 
uniform convergence of paths in an $L_2$ sense. We could take the supremum norm in $D$, since the 
approaching sequence has the same(!) jump points path wise as the limit. 

For the convergence of the discrete processes to the limit we again use a WBP on the binary 
tree. Technically we embed the Quicksort process and the discrete version all at the same time into a 
WBP on the binary tree with very specific rvs. An additional parameter $n\in \N\cup\{\infty\}$ stands 
for the input list length, $n=\infty$ for the limit. The basic idea is to find appropriate rvs on the 
binary tree such that we have a limit (of $R_n$) as $n\to \infty$. The main result of this paper is 
Theorem \ref{MainTheorem2} on weak convergence of processes. This embedding seems to be the right object
in order to obtain stronger convergence results. Notice the interplay of rvs and distributions and that 
we have to take nice versions of rvs in order to obtain (only) distributional results. We face a similar 
interplay (and technique) in the analysis of Find \cite{GruebelRoesler, RoeslerRevisited, KnoRoe}.

After this rough idea of the abstract embedding of this paper into the existing literature and of its 
contribution, we shall give some more details. 
  
Quicksort on the fly sorts the input of $n$ distinct numbers and returns first 
the smallest, than second smallest and so on at the time of identification. 
The output is a stream of data in time. We are interested in this flow of data and 
will analyze the difference to the average behavior.  

Like in Quicksort \cite{Hoare61,Hoare62}, choose with a uniform distribution a pivot from the input set $S$ 
of distinct numbers, split $S$ into the set $S_<$ of strictly smaller ones than the pivot, the pivot and the set $S_>$ 
of strictly larger ones in this order. Then continue recalling Quicksort independently always for the left list $S_<.$ 
If $S_<$ is empty continue with the next leftmost sublist recalling the algorithm. If $S_<$ consists only 
of one element, output this number immediately and continue with the next leftmost list. 
 
Let $X(S,l),l=1,2,\ldots,|S|$ denote the number of comparisons made up to the event when the $l$-th smallest 
number in the set $S$ is determined. The interpretation of $X(S,l)$ is as proportional to the 
time of publishing the $l$-th smallest. The rv $X=X(S,l)_{l=1}^{|S|}$ satisfies the recursion 
  \begin{equation}\label{1.1}
  (X(S,l))_{l=1}^{|S|} = (|S|-1 + \1_{l\le I} X(S_<,l) + \1_{l>I}(X(S_<,I-1) + X(S_>,l-I)))_{l=1}^{|S|}
  \end{equation}
for $|S|\ge 2.$ Here $I=I(S)=|S_<|+1$ denotes the rank of the pivot and has values in 
$\{1,2,\ldots,n\}$ with a uniform distribution. The rv $I=I(S)$ is independent 
of all $X$-rvs on the right of equation (\ref{1.1}). The rv $X(S_<,\cdot)$ is $0$ if $S$ is the empty set or 
has only one element. Otherwise the rv satisfies a similar recursion, where $I(S_<)$ is independent of $I(S).$ (Mathematically correct: $I(A),\ \emptyset\ne A\subset S$ are independent with $I(A)$ a uniform distribution 
on $A$ and the $X$ are recursively defined.) Continuing this way we find 
the distribution of $X(S,|S|)$ as the Quicksort distribution sorting $S$ by standard Quicksort. 

The equation (\ref{1.1}) determines the distribution of $(X(S,l))_{l=1}^{|S|}$ via the distribution for smaller   
sets. The distribution depends on $S$ only via the size $n=|S|$ of $S.$ (Prove this by induction 
on $n$ and notice the distribution of $I(S)$ depends only on the size $n$ of $S$.) For that reason we write 
$X(S,l)\stackrel{\cal D}{=}X(|S|,l).$ For notational reasons as above, which corresponds nicely to the interpretation, 
we use the boundary conditions $X(n,0)=0$ for all $n\in\N$ and $X(0,\cdot)=X(1,\cdot)\equiv 0.$ Then the above 
recursion writes for $n\ge 2$ 
  \begin{equation}\label{1.2}
  (X(n,l))_{l=1}^n \stackrel{\cal D}{=} (n-1 + \1_{l<I_n}X^1(I_n-1,l) 
    +\1_{l\ge I_n}(X^1(I_n-1,I_n-1)+X^2(n-I_n,l-I_n)))_{l=1}^n
  \end{equation}
and determines the distribution of $(X(n,l))_{l=0}^n,\ n\in\N,$ by the previous ones. The rvs 
$I_n,$ $(X^i(j,k))_{k=0}^j,\ i=1,2,\ 1\le j<n$ are independent. The rv $I_n$ has values in 
$\{1,2,\ldots,n\}$ with a uniform distribution, $X^i(j,\cdot)$ has the same distribution as 
$X(j,\cdot)$ given by recursion. Notice $X(n,n)=X(n,n-1).$ 

In our version of Quicksort we use internal randomness by picking the pivot with a uniform distribution. 
Like in standard Quicksort, we could instead of internal randomness also use external randomness. 
Choose as input an uniform distribution on all permutations $\pi$ of order $n$ and pick as pivot 
any, for example always the first in the list. Now $X(\pi,\cdot)$ is a deterministic function 
depending on the input $\pi.$ Seen as a rv with random input $\pi$ we face the same distribution as with 
internal randomness. The main advantage using internal randomness is that $X(\pi,\cdot)$ has the same distribution 
for every input $\pi$ of the same size. Alternatively we could start with $n$ iid random variables uniformly on 
$[0,1]$ as the input and choose as pivot always the first element of the list. The algorithm itself would be deterministic,
the running time is a rv via the input of an iid sequence and has the same distribution as our $X(n,\cdot).$   
 
From equation (\ref{1.2}) we obtain a recursion for the expectation $a(n,l)=E(X(n,l)),\ n\ge 2,\ 1\le l\le n$
  \begin{equation}\label{Expectation}
  a(n,l)=n-1+ \frac{1}{n}\sum_{j=1}^l (a(j-1,j-1))+a(n-j,l-j))
   + \frac{1}{n}\sum_{j=l+1}^n a(j-1,l)
  \end{equation}
The term $a(n,n)$ is the expectation of sorting $n$ numbers by Quicksort. All $a(n,l)$ are 
uniquely defined by the above equations and the starting conditions. Mart\'{i}nez \cite{Mar}
obtained the explicit formula 
  \begin{equation}\label{1.2a}
  a(n,l)=2n+2(n+1)H_n-2(n+3-l)H_{n+1-l}-6l+6
  \end{equation}
$1\le l\le n\in\N.$ $H_j$ denotes the $j$-th harmonic number $H_j=\sum_{i=1}^j\frac{1}{i}.$ 
(Notice $a(0,0)=0=a(1,0)=a(1,1)$ but in general the formula (\ref{1.2a}) is not $EX(n,l)$ 
for $n\ge 2$ and $l=0.$) For Quicksort we obtain the well known formula $a(n,n)=2(n+1)H_n-4n.$ 

Mart\'{i}nez argued with {\it Partial Quicksort} $PQ(n,l)$, which for fixed $n,l$ sorts the $l$ 
smallest elements of a list. For more results and versions of it, optimality and one-dimensional 
distributions for Partial Quicksort see \cite{MarRoe}. The Quicksort on the fly process is an extension 
of Partial Quicksort in the sense of taking $l$ as a time variable and considering processes. We 
find first up to the $l-1$-smallest elements, then continue this search for the $l$-th smallest, then 
$l+1$-smallest and so on. 

Now we come to the distribution of the process $X(n,\cdot)$ in the limit, where $X$ is defined via (\ref{1.2}).
This includes the question, how much $X(n,l)$ differs from the average $a(n,l).$ Define the rvs    
  \begin{equation}\label{1.2b}
  Y_n(\frac{l}{n})=\frac{X(n,l)-E(X(n,l))}{n}
  \end{equation}
for $l=0,1,\ldots,n$. These rvs satisfy the recursion, $n\ge 2$ 
  \begin{eqnarray*}\label{1.3}
  \left(Y_n(\frac{l}{n})\right)_{1\le l\le n}&\stackrel{\cal D}{=}&
    (C(n,l,I_n)+\1_{l<I_n}\frac{I_n-1}{n}Y^1_{I_n-1}(\frac{l}{I_n-1})\\
  &+&\1_{l=I_n}\frac{I_n-1}{n}Y^1_{I_n-1}(1) \\ 
  &+&\1_{l>I_n}(\frac{I_n-1}{n}Y^1_{I_n-1}(1) + (\frac{n-I_n}{n})Y^2_{n-I_n}(\frac{l-I_n}{n-I_n})))_{1\le l\le n} \\ \label{1.3a}
  C(n,l,i)&=& \frac{1}{n}(n-1-a(n,l)+\1_{l < i}(a(i-1,l+1)\\
  &+&\1_{l+1=i}a(i-1,i-1)+\1_{l>i}( a(i-1,i-1)+a(n-i,l-i))\\
  C(n,0,i)&=&0  
  \end{eqnarray*}
for $l=1,\ldots,n-1.$ Notice $Y_n$ is well defined, with the help of the indicator function, and 
there are no boundary conditions besides $Y_0\equiv 0\equiv Y_1.$  

We extend the process $Y_n$ nicely to a process on the unit interval $[0,1]$ with values in the 
space $D=D[0,1]$ of cadlag functions (right continuous functions with existing 
left limits \cite{Billingsley}) on the unit interval. This can be done by linear interpolation or a 
piece wise constant function. We shall use the extension
  \begin{equation}\label{1.3b}
  Y_n(t):=Y_n(\frac{\lfloor n t\rfloor}{n}).
  \end{equation}
The process $Y_n$ with values in $D$ satisfies the recursion, we use $U_n=\frac{I_n}{n}$ 
  \begin{eqnarray*}
  Y_n&\stackrel{\cal D}{=}& ( C(n,\lfloor n t\rfloor,I_n)+\1_{t < U_n} \frac{I_n-1}{n}Y^1_{I_n-1}(\frac{n t}{I_n-1}\wedge 1)\\
  && + \1_{t \ge U_n}(\frac{I_n-1}{n}Y^1_{I_n-1}(1) + (\frac{n-I_n}{n}) Y^2_{n-I_n}(\frac{t - U_n}{1-U_n})))_{t\in [0,1]}
  \end{eqnarray*}
for $n\in\N.$ In short notation 
  \begin{equation}  
  Y_n\stackrel{\cal D}{=}\varphi_n(U_n,(Y^1_k)_{k<n},(Y^2_k)_{k<n})
  \end{equation}
for a suitable function $\varphi_n.$  

If $n\to\infty$ then $U_n$ converges in distribution to a rv $U$ with a uniform distribution. 
We might expect that the process $Y_n$ converges in some sense to a limiting 
process $Y=(Y(t))_{t\in [0,1]}$ with values in $D$ satisfying something like 
the stochastic fixed point equation
  \begin{eqnarray}\nonumber
  Y&\stackrel{\cal D}{=}& (\1_{t< U}U Y^1(\frac{t}{U}) + \1_{t\ge U}(UY^1(1)+(1-U)Y^2(\frac{t-U}{1-U}))
    +C(t,U))_t\\
  Y&\stackrel{\cal D}{=}&\varphi(U,Y^1,Y^2)\label{1.4}
  \end{eqnarray}
for a suitable function $\varphi.$ The rvs $Y^1,Y^2,U$ are independent. $Y^1$ and $Y^2$ have 
the same distribution as $Y$ and $U$ is uniformly distributed on the unit interval $[0,1].$ 
The cost function $C=C(\cdot,U)$ is given by
  \begin{eqnarray}\label{1.5}
  C(t,x)&:=& C(x)+2\1_{t< x}(-1+x+(1-t)\ln(1-t)-(1-x)\ln(1-x)\nonumber\\
      &&-(x-t)\ln(x-t))\\
  C(x)&:=& 1+2x\ln x+2(1-x)\ln(1-x)
  \end{eqnarray}
and is the limit of $C(n,l_n,i_n)$ as $n\to\infty$ with $\frac{l_n}{n}\to_n t,\frac{i_n}{n}\to_n x,$ Proposition
(\ref{5.5}).

Our first major result, Theorem \ref{Hauptsatz}, states the existence of the $Y$-process with values in $D$. 
  \begin{Theorem}\label{MainTheorem}
  Let $U^v,\ v\in V,$ be iid rvs with a uniform distribution on $[0,1]$ and $V$ be the binary tree.
  Then there exists a family $Y^v,\ v\in V,$ of rvs with values in $D$ and all of the same 
  distribution satisfying almost surely 
    \begin{eqnarray}\label{1.6}
    Y^v&=& \varphi(U^v,Y^{v1},Y^{v2})\\ \nonumber
    Y^v(t)&=&\1_{t< U^v}U^v Y^{v1}(\frac{t}{U^v}) + \1_{t\ge U^v}(U^v Y^{v1}(1)+(1-U^v) Y^{v2}(\frac{t-U^v}{1-U^v}))+C(t,U^v)
    \end{eqnarray}
  simultaneously for all $t\in[0,1].$ 
  \end{Theorem}
The above equation (\ref{1.6}) for rvs implies the distributional equality and $Y=Y^\emptyset$ satisfies the 
fixed point equation (\ref{1.4}). A specific family $Y^v$ is explicitly given in the paper. 
We call $Y^\emptyset$ the Quicksort process and its distribution the Quicksort process distribution.
 
Our second major result, compare Theorem \ref{4.8} and Corollary \ref{finite dimensional Convergence}, states the 
weak convergence of $Y_n$ to the Quicksort process constructed above.

  \begin{Theorem} \label{MainTheorem2}
  All finite dimensional distributions of $Y_n$ without the coordinate zero converge to those of 
  the Quicksort process $Y$ constructed above. 
  \end{Theorem}

The first result is a probabilistic result, while the second is a measure theoretic one. We 
obtain both results via the Weighted Branching Process \cite{RoeslerBielefeld} and an explicitly given 
nice family of processes $Y_n^v$ indexed by $n\in\N$ and the binary tree. Basically we 
use the splitting $U$-rvs for the $Y$-process also for the $Y_n$-process. 

Inspired by previous work on Find \cite{KnoRoe} one might choose the pure measure theoretic approach by 
first showing the convergence of finite dimensional distributions of $Y_n$ by the contraction method obtaining 
a limiting measure on $\R^{[0,1]}.$ Then show the tightness \cite{Billingsley} of the sequence of measures in 
order to obtain a limiting measure with outer measure $1$ of $D.$ The measure restricted to $D$ provides 
the desired distribution. There exist rvs satisfying (\ref{1.4}). (Statement (\ref{1.6}) is slightly stronger.) 
Neininger and Sulzbach \cite{NeiSul} study this approach in more generality using the Zolotarev metric. 

We prefer here a probabilistic approach using rvs instead of measures. Actually we formulated a (seemingly) 
stronger statement, the existence of a family $Y^v$ of rvs with values in $D$ and indexed by the binary 
tree, which satisfies (\ref{1.4}) almost everywhere as rvs. For the construction we used a.e. convergence 
of appropriate rvs in supremum metric. A similar probabilistic approach via the Skorodhod metric was used 
\cite{GruebelRoesler} for Find, the first process analysis of a stochastic algorithm. 

Also for Quicksort we find the two approaches in the literature with different types of results. The 
probabilistic approach by R\'egnier \cite{Regnier} provided a limiting rv via an $L_2$-martingale. The 
measure theoretic approach via the contraction method \cite{Roesler92} and the backward view characterized 
the limiting distribution as unique solution of a stochastic fixed point equation. Via the approach 
of Weighted Branching Processes one can construct a family of tree indexed rvs satisfying the stochastic 
fixed point equation as rvs. (Another $L_2$-martingale $R_n=\sum_i L_v C^v$, different from R\'egnier's, is 
the key to convergence.)

We stated the results for processes with values in $D.$ This is because $D$ is preferred 
to the space $E$ of left continuous functions $f:[0,1]\to\R$ with existing right limits. 
Both spaces are isomorphic with respect to the topological, probabilistic and algebraic 
structure. (The easiest argument is via the map $[0,1]\ni x\mapsto 1-x.$) The space $D$ 
is preferred in probability theory (\cite{Billingsley}) and this is our main argument for using 
the space of cadlag functions. This explains the motivation for our normalization in (\ref{1.2b}).

Using the same arguments, the path wise left continuous version of $Y$ is a solution in $E$ 
of the corresponding equation (\ref{1.6}) in $E.$ For the discrete setting and interpretation 
it seems on the first view more natural to use left continuous functions. If we know 
the $l$-th largest and are only interested in that one, we do not have to look to the future, 
we are done. The $D$ and $E$ versions are mathematically equivalent.

\section{The Weighted Branching Process}
We introduce Weighted Branching Processes (WBP) in general and will specialize to Quicksort examples
on the binary tree. 

Let $V$ be the Ulam-Harris tree ${\N}^*=\cup_{n=0}^\infty{\N}^n$ of all finite sequences of 
natural numbers including the empty sequence denoted by the empty set $\emptyset.$ (By convention
${\N}^0=\{\emptyset\}.$) $V$ is a rooted tree with root $\emptyset$ and the edges 
$(v,vi),\ v\in V,\ i\in \N$ in graph theoretical sense. We use the standard notation 
$v=(v_1,v_2,v_3,\ldots,v_n)=v_1v_2\ldots v_n$ for a vertex $v\in V$ of length 
$n=|v|.$ We skip the empty set in the notation whenever possible. We use $V_n$ for $v\in V$ of length $n$ 
and $V_{\le n},V_{<n}$ appropriate. 

A weighted branching process (WBP) is a tuple $(V,(G,*,H,\odot),(T,C)).$ $(T,C)$ is a random variable 
with values in $G^{\N}\times H$ on some probability space. (Actually we need only this distribution
for the description, but prefer here a probabilistic language.) $(G,*)$ is a measurable semi group 
($*:G\times G\to G,\ (g,h)\mapsto g*h$ associative and measurable) with a neutral element $e$ 
($\forall g\in G:\ e*g=g=g*e$) and a grave $\triangle$ ($\forall g\in G:\ \triangle*g=\triangle=g*\triangle,$
once in the grave forever in the grave). $(G,*)$ operates left on $H$ via $\odot:G\times H\to H.$ 
$H$ is a measurable space and $G\times H$ endowed with the product $\sigma$-field. 

Let $(T^v,C^v),\ v\in V,$ be independent copies of $(T,C)$ on the same probability space $(\Omega,\A,P).$
The interpretation of $C^v$ is as a cost function on a vertex $v$ and $T^v_i,$ where $T=(T_1,T_2,\ldots),$ 
is a transformation (weight) on the edge $(v,vi).$ The interpretation of $G$ is as a map from $H$ to $H$. 
If $H$ has an additional structure then we might enlarge $G$ to have the induced structure. Examples are 
$H$ is a vector space or an ordered set and the extended $G$ will be a vector space or ordered set 
via the natural extension.   

For a WBP define the path weights $L_w^v:\Omega \to G,\ v,w\in V$ on paths $(v,vw)$ recursively 
by $L^v_\emptyset=e$ and 
  $$L^v_{w i}=L^v_w*T_i^{v w}.$$

One of the basic assumptions of a WBP is the independence of families, but arbitrary dependence within a family.
Let $\A_n$ be the $\sigma$-field generated by all $T_i^v,C^v,$ for $v\in V_{<n},i\in\N.$ Then $L^v_w$ is 
measurable with respect to $\A_{|v w|}$ and $C^v$ is independent of $\A_{|v|}.$ We will use this many times
in the sequel.  

For a WBP without costs we write also $(V,(G,*),T).$ We drop the vertex $\emptyset$ whenever possible, 
e.g. $L_v=L_v^\emptyset$ or later $R^\emptyset=R,\ S^\emptyset=S.$ The interpretation of $L^v_w$ the grave is, 
we can not see the path from the root $v$ to $vw.$ Mathematically, no value grave connected to the path $(v,vw)$ 
will ever contribute, like $L^v_w C^{vw}$ will be $0$ in our examples if $L^v_w=0.$ By this construction we shall use 
freely other trees like $m$-ary trees $\{1,2,\ldots,m\}^*=\cup_{n\in\N_0}\{1,2,\ldots,m\}^n$ of all finite sequences 
over the finite alphabet $1,2,\ldots,m$ in an appropriate sense. For the $m$-ary tree we take $T=(T_1,T_2,\ldots,T_m)$ instead of $T=(T_1,T_2,\ldots,T_m,\triangle,\triangle,\ldots)$ on the original tree $\N^*.$ 

For the next sections we need the following two examples. Although they provide known results 
the novelty is the line of arguments, which can be generalized and which are the key for the Quicksort process. 
 
{\bf Example 1:} Quicksort distribution \cite{Roesler91}: Here we show mainly the existence of 
the Quicksort distribution via an embedding into the WBP. Consider the WBP $(V=\{1,2\}^*,(\R,\cdot,\R,\cdot),((U,1-U),C(U)))$ 
with $U$ has a uniform distribution and   
  $$
  C(x)=1+2x\ln x+2(1-x)\ln(1-x)
  $$
as in (\ref{1.5}).

$G$ is the multiplicative semi group $\R$ with the neutral element $e=1$ and 
the grave $\triangle=0.$ $G$ operates left on $H=\R$ by multiplication. 
Let $U^v,\ v\in V,$ be independent rvs with a uniform distribution on $[0,1]$. Put
  $$T_1^v=U^v,\quad T_2^v=1-U^v,\quad  C^v=C(U^v).$$
(For the general WBP with tree $\N^*$ we could take $T_i^v=\triangle=0$ for $i\ge 3$. The 
smaller binary tree is more suitable here.) Since $H$ is an ordered vector space, we 
extend $G$ with the interpretation of maps to the ordered vector space generated by the maps. 

The total weighted cost $R_m:=\sum_{v\in V_{<m}}L_v C^v$ up to the $m-1$ generation is an 
$L_2$-martingale and converges in $L_2$ and a.e. to a rv $Q$ \cite{Roesler91}. The distribution 
of $Q$ is called the Quicksort distribution. The distribution is uniquely characterized 
\cite{Roesler91,FillJanson} as the solution of the stochastic fixed point equation 
  \begin{equation}\label{2.2a}
  Q\stackrel{\cal D}{=}UQ_1+(1-U)Q_2+C(U)
  \end{equation}
with expectation $0.$ Here $\stackrel{\cal D}{=}$ denotes equality in distribution.
The random variables $U,Q_1,Q_2$ are independent, $U$ is uniformly distributed and $Q_1,Q_2$ 
have the same distribution as $Q.$ 

By the a.s. convergence of $R_m^v=\sum_{w\in V_{<m}}L^v_w C^{v w}$ the rvs 
  \begin{equation}\label{2.2}
  Q^v:=\sum_{w\in V}L^v_w C^{v w}
  \end{equation}
exist and satisfy a.e. 
  \begin{equation}\label{2.3}
  Q^v=U^v Q^{v1}+(1-U^v)Q^{v2}+C(U^v)
  \end{equation} 
for every $v\in V.$ Of course the distribution of $Q^v$ is a solution of (\ref{2.2a}) and is the 
Quicksort distribution.

{\bf Example 2:} Convergence of the discrete Quicksort distributions \cite{Roesler91}: The original 
problem concerns the number $X_n$ of comparisons to sort $n$ distinct reals. We use 
internal randomness. Then for $n\in\N$ 
  $$X_n\stackrel{\cal D}{=}n-1+X^1_{I_n-1}+X_{n-I_n}^2$$
with $I_n,X^1,X^2$ are independent, $I_n$ has a uniform distribution on $1,\ldots,n$ 
and $X^1_i,X^2_i$ have the same distribution as $X_i.$ The boundary conditions are 
$X_0$ and $X_1$ are identical $0.$ The expectation of $X_n$ is $a_n=a(n,n)$ as in (\ref{1.2a}).
The normalized rvs $Y_n=\frac{X_n-a_n}{n}$ satisfy the recursion
  $$Y_n\stackrel{\cal D}{=}\frac{I_n-1}{n}Y^1_{I_n-1}+\frac{n-I_n}{n}Y^2_{n-I_n}+C_n(I_n)$$
where 
  $$C_n(i):=\frac{n-1-a_n+a_{i-1}+a_{n-i}}{n}$$

We come now to the abstract embedding of this example into a WBP with an additional parameter $n\in\N.$ 
Let $H$ be the set of functions $h:\N_0\to\R$ and $G$ the set $H\times G_2$ where 
$G_2$ are the functions $g:\N_0\to\N_0$ satisfying $g(0)=0$ and $g(n)<n$ for all 
$n\in\N.$ The semi group structure is given by 
  $$(f_1,g_1)*(f_2,g_2)=(f_1f_2\circ g_2, g_2\circ g_1)$$
and the operation on $H$ via 
  $$(f,g)\odot h=f h\circ g\ \ \ \ ((f,g)\odot h)(n)=f(n)h(g(n))$$ 
($\circ$ denotes composition and we use multiplication on $\R$).

The interpretation of $(f,g)\in G$ is as a map on $H,$ where $f$ is a multiplicative factor and
$g$ an index transformation. The operation $*$ corresponds to the convolution  
of maps on $H.$ Since $H$ is a vector space we may enlarge $G$ naturally to a vector space. 

Consider the binary tree $V=\{1,2\}^*$ and let $U^v,v\in V,$ be independent rvs with a uniform distribution. 
Let $I^v_n:=\lceil n U^v\rceil$ (upper Gauss bracket) and define the transformations on the 
edges $(v,v1),(v,v2)$ by 
  \begin{eqnarray*}
  J^v_1(n):=I_n^v-1  &\quad & J_2^v(n):=n-I_n^v\\
  T_1^v(n):=(\frac{J^v_1(n)}{n},J^v_1(n))  
    &\quad& T_2^v(n):=(\frac{J_2^v(n)}{n},J_2^v(n))
  \end{eqnarray*}
and the vertex weight 
  $$C^v(n):=C_n(I_n^v)$$
The rvs 
  $$R_m^v:=\sum_{w\in V_{<m}}L^v_w\odot C^{v w}$$
converge as $m\to \infty$ a.e. and in $L_2$ to a limit $R^v$ and satisfy 
  $$R_m^v=\sum_i T_i^v \odot R_{m-1}^{vi}+C^v$$
for $m\in\overline {N}=\N\cup\{\infty\}.$ 

Notice $R_m^v$ and $R^v$ take values in $H$ and $R^v(n),\ n\in\N,$ 
is a random variable with values in the reals. Notice the connection to the previous description
of the Quicksort rv. $Y_n$ from the introduction   
  $$Y_n\stackrel{\cal D}{=}R^\emptyset(n).$$
$R^v(n)$ converge for every $v\in V$ in $L_2$ to the rv $Q^v$ from the Quicksort 
example \cite{Roesler91}. We shall use $Q^v_n=R^v(n)$ in the sequel and drop the root whenever possible. 

{\bf Example 3:} Joint embedding: It is worth while to put the two examples together. Use 
$\overline{\N_0}=\N_0\cup \{\infty\}$ instead of $\N_0$ in the second example and incorporate the 
first example via the value $\infty.$

\section{The Quicksort Process}
In this section we consider the Quicksort process. Let $D=D([0,1])$ be the vector space of 
cadlag functions $f:[0,1]\to\R$ (right continuous with existing left limits).
$D$ is endowed with the Skorodhod topology \cite{Billingsley} induced by the 
Skorodhod $J_1$-metric 
  \begin{equation}\label{3.1}
  d(f,g)=\inf\{\epsilon>0\mid \exists \lambda\in \Lambda: \|f-g\circ \lambda\|_\infty<\epsilon,\   
    \|\lambda-\mbox{id}\|_\infty<\epsilon\}
  \end{equation}
where $\Lambda$ is the set of all bijective increasing functions $\lambda:[0,1]\to [0,1]$.
We use the supremum norm $\|f\|_\infty=\sup_t|f(t)|.$  

The space $(D,d)$ is a separable, non complete metric space, but a polish space \cite{Billingsley}. 
The $\sigma$-field $\sigma(D)$ is the Borel-$\sigma$-field via the Skorodhod metric. 
The $\sigma$-field is isomorphic to the product $\sigma$-field $\R^A\cap D$ where 
$A$ is a dense subset of $[0,1]$ containing the $1.$ 

Let $\F(D)$ be the space of all measurable functions $X$ with values in $D.$ For $1\le p<\infty$ 
let $\F_p(D)$ be the subspace such that 
  \begin{equation}\label{3.2}
  \|X\|_{\infty,p}:=\| \|X\|_\infty \|_p<\infty
  \end{equation}
is finite. Here $\|\cdot\|_p, p<\infty,$  denotes the usual $L_p$-norm for rvs. The map 
$\|\cdot\|_{\infty,p}$ is a pseudo metric on $\F_p(D).$ Let $\sim$ be the common equivalence relation
  $$X\sim Y\Leftrightarrow P(X\ne Y)=0$$
and $F_p(D)$ be the set of equivalence classes $[X]=\{Y\in \F(D)\mid X\sim Y\}$ intersected with 
$\F_p(D).$ Then it is well known
  \begin{Proposition}
  $(F_p(D),\|\cdot\|_{\infty,p})$ is a Banach space for $1\le p<\infty$ with the usual addition 
  and multiplication
    $$[f]+[g]=[f+g],\quad c[f]=[cf],\quad \|[f]\|_{p,\infty}=\|f\|_{\infty,p}.$$
  \end{Proposition} 
In the following we will be careless and will not differ between functions and equivalence classes. 

Let $(G,*)$ be the semi group $G=D\times D_{\uparrow}$ where $D_\uparrow$ consists of increasing 
functions $D\ni g:[0,1]\to[0,1]$ and the semi group operation $*$ is 
  $$(f_1,g_1)*(f_2,g_2):= (f_1 f_2\circ g_1,g_2\circ g_1).$$
$(G,*)$ has as neutral element $(1,\mbox{id})$, the function identically $1$ and the identity, 
and the grave is $(0,\mbox{id}).$ $G$ operates left on $H=D$ via 
  $$(f,g)\odot h:=f h\circ g.$$
The tuple $(f,g)\in G$ has the interpretation of a map $M_{f,g}$ from $H$ to $H$ acting as 
  $$M_{f,g}(h)(t)=f(t)h(g(t)).$$
The first coordinate $f$ acts as a space transformation, the second coordinate $g$ as a time 
transformation. The semi group structure $*$ is the composition of the corresponding maps 
  $$M_{(f_1,g_1)*(f_2,g_2)}=M_{f_1,g_1}\circ M_{f_2,g_2}.$$
(Notice the order of the composition.) Since $H$ is a vector space and $\R$ is a lattice, we 
will embed $G$ to maps $H^H$ and use freely the induced structure $+,\cdot,\vee$ on $G,$ $\vee$ 
denotes the supremum, $a\in\R$  
  $$(M_{f_1,g_1}+M_{f_2,g_2})(h)=M_{f_1,g_1}(h)+M_{f_2,g_2}(h)$$
  $$ a\cdot  (M_{f,g}(h))= (a\cdot M_{f,g})(h)$$
  $$(M_{f_1,g_1}\vee M_{f_2,g_2})(h)= (M_{f_1,g_1}(h))\vee (M_{f_2,g_2}(h)).$$

Let $V=\{1,2\}^*$ be the binary tree and let $(U^v,Q^v),\ v\in V,$ be the rvs as in the Quicksort 
example in the WBP section. Define on the edges $(v,vi),\ v\in V,i\in\{1,2\}$ the edge weights 
$T^v=(T_1^v,T_2^v),\ T_i^v=(A_i^v,B_i^v)$ with values in $G$ and the vertex weights $C^v$ 
with values in $H$ by 
  $$A_1^v=\1_{[0,U^v)}U^v \qquad  B_1^v(t)=1\wedge\frac{t}{U^v}$$
  $$A_2^v=\1_{[U^v,1]}(1-U^v)\qquad  B_2^v(t)=0\vee \frac{t-U^v}{1-U^v}$$
  $$C^v=C(\cdot,U^v)+\1_{[U^v,1]}Q^{v1}$$  
where $C(t,x)$ is given in (\ref{1.5}). Define rv $R_m^v$ with values in $H$ by 
  \begin{equation}\label{Rm} 
  R^v_m:=\sum_{w\in V_{<m}}L_w^v\odot C^{v w}
  \end{equation}
for $m\in\N_0.$
  \begin{Proposition}\label{3.5}
  The rvs $R_m^v$ satisfy  
    \begin{equation}\label{3.6}
    R_m^v=\sum_{i=1}^2T_i^v R_{m-1}^{vi}+C^v
    \end{equation}
for all $m\in\N,\ v\in V.$ 
  \end{Proposition}
The proof is straight forward and easy since $*$ and $+$ interchange. \hfill q.e.d.   
   
  \begin{Theorem} \label{Hauptsatz} Let $(V=\{1,2\}^*,(G,*,H,\odot),((T_1,T_2),C))$ be the weighted 
  branching process as defined above. Then $R_m^v$ converges in supremum metric on $D$ as $m\to\infty$ a.e. to a 
  rv $R^v$ for all $v\in V.$ The family $R^v,\ v\in V $ satisfies 
  \begin{equation}\label{Fixpunkt}
    R^v=\sum_{i=1}^2T^v_i R^{vi}+C^v
  \end{equation}
  almost everywhere. For every $p>1$ and $v\in V$ holds 
    \begin{equation}\label{Maximum}
    \|\|R_m\|_\infty\|_p\le \frac{8+2^{-1/p}k_p\|Q\|_p}{1-k_p}
    \end{equation}
  where $k_p=\left(\frac{2}{p+1}\right)^{1/p}$ and $Q$ is a rv with the Quicksort distribution.
  \end{Theorem}
Proof: Let $S_m^v=R_{m+1}^v-R^v_m=\sum_{w\in V_m} L_w^v\odot C^{v w},\ S_0^v=C^v$ and $b_m:=\|S_m^\emptyset\|_{2,\infty}.$ Notice 
$b_m=\|S_m^v\|_{\infty,2}$ does not depend on the vertex $v$ and 
  $$S_m^v=\sum_{i=1}^2 T_i^v*S^{vi}_{m-1}$$ 
by (\ref{3.6}) for all $m\in \N,\ v\in V.$ 

$\bullet$ $b^2_m\le \frac{2}{3}b^2_{m-1}$ for $m\in\N.$ 

  \begin{eqnarray*}
  b_m^2&=&E\sup_t(\1_{t< U}US_{m-1}^1(\frac{t}{U})+\1_{t\ge U}(1-U)S_{m-1}^2(\frac{t-U}{1-U}))^2\\
  &=&E\sup_t((\1_{t<U}US_{m-1}^1(\frac{t}{U}))^2+(\1_{t\ge U}(1-U)S_{m-1}^2(\frac{t-U}{1-U}))^2)\\
  &\le& E((U\|S_{m-1}^1\|_\infty)^2)+E((1-U)\|S_{m-1}^2\|_\infty)^2)\\
  &=& (E(U^2)+E((1-U)^2))b_{m-1}^2=\frac{2}{3}b_{m-1}^2
  \end{eqnarray*}
     
$\bullet$ $b_m\le \left(\frac{2}{3}\right)^{m/2}b_0$ where $b_0=\|C\|_{\infty,2}<\infty.$ 
    
Easy by recursion. 

$\bullet$ $R_m,\ m\in\N,$ is a Cauchy sequence in $(\F_2(H),\|\cdot\|_{\infty,2}).$ 

For $n_0\le m<n$ argue     
  \begin{eqnarray*}
  \|R_n-R_m\|_{\infty,2}&\le& \|\sum_{l=m}^{n-1}S_l\|_{\infty,2}\\
  &\le& \sum_{l=m}^{n-1}\|S_l\|_{\infty,2}\le \sum_{l\ge n_0}b_l\\
  &\le &b_0\left(\frac{2}{3}\right)^{n_0/2}\frac{1}{1-\sqrt{2/3}}\to_{n_0\to\infty}0 
  \end{eqnarray*}
  
$\bullet$ $\sum_m\|S_m\|_\infty<\infty$ a.e.

  $$E(\sum_m\|S_m\|_\infty)=\sum_m \|\|S_m\|_\infty \|_1\le \sum_m\|\|S_m\|_\infty\|_2<\infty$$

Let $R$ be the limit of the Cauchy sequence $R_m$. $R(t)$ is the point wise limit
of $\sum_{m\in\N}S_m(t)$ for all $t\in[0,1].$ 

$\bullet$ $R$ is well defined a.e..

Easy by the previous statements. 

$\bullet$ $R\in\F_2(D).$ 

Every $R_m$ is in $\F_2(D).$ Since we can write $R=\sum_{m=0}^\infty S_m$ a.e. 
the triangle inequality provides the result. 

$\bullet$ Equation (\ref{Fixpunkt}) is true.  

The same is true for every $R_m^v,\ v\in V,$ instead of $R_m$ and call $R^v$ the limit in a.s. 
sense. Then $R_m^v$ satisfies (\ref{3.6}) and the a.e. convergence provides (\ref{Fixpunkt}).  

For the last statement notice 
  $$\|R_m\|_\infty\le \|(T_1R_{m-1}^1)\vee (T_2R_{m-1}^2)\|_\infty+\|C\|_\infty.$$
Consider the WBP as above but with cost function the constant $\overline C^v:=8+U^v|Q^v|,\ v\in V,$ instead of $C^v.$ Define 
  $$\overline R^v_m:=\sum_{j<m}\bigvee_{w\in V_j}L^v_w\overline C^{v w},$$  
where the symbol $\bigvee$ denotes the supremum. Then $\overline R^v_0=0$ and $\overline R^v_m$
increases in $n$ point wise to $\overline R^v$ for all $v\in V$ and 
  $$\overline R_m^v=T_1^v\overline R_{m-1}^{v1}\vee T_2^v\overline R_{m-1}^{v2}+\overline C^v$$
for $m\in\N.$ By induction it is easy to show $\|R_m^v\|_\infty\le \overline R^v_m.$
(Notice $\|C\|_\infty\le 8.$) Then for $p>1$        
  \begin{eqnarray*}
  \|\overline R_m\|_p&\le& \|T_1\overline R_{m-1}^1\vee T_2\overline R_{m-1}^2\|_p+\|\overline C\|_p\\
  &\le& \left(E|T_1\overline R_{m-1}^1|^p+E|T_2\overline R_{m-1}^2|^p\right)^{1/p}+ 8+\|UQ\|_p\\
  &\le& \|\overline R_{m-1}\|_p\left(E|T_1|^p+E|T_2|^p\right)^{1/p}+ 8+\|UQ\|_p\\
  &\le&  k_p\|\overline R_{m-1}\|_p+ 8+\|UQ\|_p\\         
  &\le& \sum_{i=0}^m k_p^i(8+\|U\|_p\|Q\|_p)\le \frac{8+\|U\|_p\|Q\|_p}{1-k_p}\\
  \|\overline R\|_p&=&\lim_m\|\overline R_m\|_p\le\frac{8+(\frac{1}{p+1})^{1/p}\|Q\|_p}{1-k_p}
  \end{eqnarray*}
\hfill q.e.d.

\section{Convergence of the discrete Quicksort process.}
In this section we prove the convergence of finite dimensional marginals of $Y_n$ to $Y.$ 
We will define a nice version of $Y_n$ such that in $L_2$-norm $Y_n(t)$ converges to $Y(t)$ 
for every $t\in[0,1].$ This requires to define a nice family $(Y_n^v)_n$ of random processes 
with values in $D,$ indexed by the tree $V$. We will include the Quicksort process via the 
index $\infty$ and consider $(Y_n)_{n\in\overline N_0},\ \overline N_0=\N\cup \{0,\infty\}.$ 
Compare this construction to the examples 3 of the section on the WBP. We shall use the 
general notation of a WBP with binary tree, but now on a more general function space. 

Let $V=\{1,2\}^*$ be the binary tree and   
  $$H=\{h:[0,1]\times\overline{\N}_0\to\R\mid \forall n\in \overline{\N}_0:\ h(\cdot,n)\in D\}.$$
Let $G_2$ be the set of all $g:[0,1]\times \overline{\N_0}\to [0,1]\times \overline{\N_0}$
such that $g(\cdot,0)\equiv 0,\ \forall n\in\N:\ \Phi_2(g(\cdot,n))<n$ and 
$\Phi_2(g(\cdot,\infty))=\infty$ where $\Phi_2$ denotes the projection to the second coordinate.   

Define $G=H\times G_2$ with the semi group operation $*$ 
  $$(f_1,g_1)*(f_2,g_2):= (f_1f_2\circ g_1,g_2\circ g_1).$$
$(G,*)$ has the neutral element $(1,\mbox{id})$, the function identically $1$ and the identity, 
and the grave is $(0,\mbox{id}).$ $G$ operates left on $H$ via 
  $$(f,g)\odot h:=f h\circ g.$$
The tuple $(f,g)$ has the interpretation of a map $M_{f,g}$ from $H$ to $H$ via 
$(M_{f,g}(h))(t,n)=f(t,n)h(g(t,n))$. The first coordinate $f$ acts as a space transformation, 
the second coordinate $g$ as a time and index transformation. Since $H$ is a vector space and ordered 
set we will embed $G$ into $H^H$ and use freely the operations $+,\cdot, \vee$ 
  $$(M_{f,g}+M_{f_1,g_1})(h)=M_{f,g}(h)+M_{f_1,g_1}(h)$$
  $$a\cdot  (M_{f,g}(h))= (a\cdot M_{f,g})(h)$$
  $$(M_{f,g}\vee M_{f_1,g_1})(h)= ((M_{f,g}(h))\vee ((M_{f_1,g_1})(h))$$
for $a\in \R.$ 
 
Let $U^v,\ v\in V,$ be independent rvs with a 
uniform distribution. Let $Q^v,Q_n^v,v\in V,$ be the rvs as in the Quicksort examples 1 and 2
in the WBP section. Define on the edges $(v,vi),v\in V,i\in\{1,2\},$ the edge weights 
(transformations) $T^v=(T_1^v,T_2^v),\ T_i^v=(A_i^v,(B_i^v,J_i^v))$ with values in $G$ and the vertex weights $C^v$ 
with values in $H$ by
  \begin{eqnarray*}
  I_n^v:=\lceil n U^v\rceil &\quad & U_n^v:=\frac{I_n^v}{n} \\
  J_1^v(t,n):=I_n^v-1&\quad & J_2^v(t,n):=n-I_n^v\\
  A_1^{v}(t,n):=\1_{t<U_n^v}\frac{I_n^v-1}{n} &\qquad &  A_2^v(t,n):=\1_{t\ge U^v_n}(1-\frac{I_n^v}{n})\\
  B_1^v(t,n):=1\wedge \frac{\lfloor n t \rfloor}{I^v_n-1}&\quad &
    B_2^v(t,n):=0\vee \frac{t-U_n^v}{1-U^v_n}\\  
  C^v(t,n):=C(n,\lfloor n t \rfloor,I_n^v)&+&\1_{t\ge U_n^v}\frac{I_n^v-1}{n}Q^{v1}_{I_n^v-1}   
  \end{eqnarray*}
for $n\in\N$ and 
  \begin{eqnarray*}
  J_1^v(t,\infty):=\infty&\quad & J_2^v(t,\infty):=\infty\\
  A_1^{v}(t,\infty):=\1_{t<U^v}U^v &\qquad &  A_2^v(t,\infty):=\1_{t\ge U^v}(1-U^v)\\
  B_1^v(t,\infty):=1\wedge \frac{t}{U^v}&\quad &
    B_2^{v}(t,\infty):=0\vee \frac{t-U^v}{1-U^v}\\  
  C^v(t,\infty):=C(t,U^v)&+&\1_{t\ge U^v}U^v Q^{v1}   
  \end{eqnarray*}
$t\in[0,1].$ $A_i^v(t,0)$ is identically $0$.

Define the rvs $R_m^v$ 
  $$R^v_m:=\sum_{w\in V_{<m}}L_w^v\odot C^{v w}$$
for $m\in\overline{\N}_0,\ v\in V.$ Notice, $R^v_m$ takes values in $H$ and 
$R^v_m(\cdot,\infty)$ is the same as the previous Quicksort $R^v_m$ given in (\ref{Rm}).
The embedding analogous to the embedding in example 3 of the WBP section. 
  \begin{Proposition}\label{4.5}
  The rvs $R_m^v$ satisfy  
    \begin{equation}\label{4.6}
    R_m^v=\sum_{i=1}^2T^v_i\odot R_{m-1}^{vi}+C^v
    \end{equation}
for all $m\in\overline{\N},\ v\in V.$ 
  \end{Proposition}
(By convention $\infty-1=\infty.$) The proof is easy since $*$ and $+$ interchange. \hfill q.e.d. 
   
For every $n\in\N$ the function $R_m^v(\cdot,n)$ converges a.e. as $m\to\infty$ to 
$R^v(\cdot,n)$ for every $v\in V.$ Notice the number of summands increase in $m$ and 
$R^v(\cdot,n)$ has only finitely many non zero summands. By induction it is easy to show 
  \begin{Proposition}\label{4.7}
    $$Y_n^v=R^v(\cdot,n)=\sum_{w\in V}L_w^v\odot C^{v w}(\cdot,n)$$
  for $n\in\N_0,v\in V.$ 
  \end{Proposition}
Proof: The cases $n=0$ and $n=1$ are easy, since both sides are $0.$  For the induction step use 
the representation after (\ref{1.3b}) for $Y_n^v$ and (\ref{4.6}) for $R^v.$ We show 
as an example the equality for the first term 
  \begin{eqnarray*}
  T_1^v\odot R^{v}(t,n)&=&A_1^v(t,n)R^{v1}(B_1^v(t,n),J_1^v(t,n))\\
  &=&A_1^v(t,n)Y^{v1}_{J_1^v(t,n)}(B_1^v(t,n))\\
  &=&\1_{t<U_n^v}\frac{I_n^v-1}{n}Y^{v1}_{I_n^v-1}(1\wedge\frac{\lfloor n t\rfloor}{I_n^v-1})
  \end{eqnarray*}
The rest follows the same line. \hfill q.e.d.  

  \begin{Theorem}\label{4.8}
  In the above setting
    $$\|Y_n^v(t)-Y^v(t)\|_2\to_{n\to \infty}0$$
  for all $t\in (0,1]$ and $v\in V.$ 
  \end{Theorem} 
Proof: Let $S_m^v=R_{m+1}^v-R_m^v$ and $a:=\sup_{t\in[0,1]}\sup_{n\in\N_0}E((C(t,n))^2).$
In the following let $v,w,\overline w\in V,\ m,\overline n\in\overline{\N}_0, t\in[0,1].$ 

Notice $L^v_w=(A_w^v,B_w^v,J_w^v)$ acts as a map on $H$ via $(L^v_w\odot h)(t,n)=:A^v_w(t,n)h(B^v_w(t,n),J^v_w(t,n)).$
(We take this as the definition of $(A_w^v,B_w^v,J_w^v).$) We will use $E(L^v_w)$ also as an operator acting on $H$
in the sense $E((L^v_w\odot h)(t,n)=((E(L^v_w))(h))(t,n).$
We use $(L)^2$ of an operator via $((L)^2)(h)=(L(h))^2.$
 
Let $\A_m$ be the $\sigma$-field generated by all $U^v,\ v\in V_{<m}$. The rv $L^v_w$ is 
measurable with respect to $\A_{|v w|}$ and $C^{v}$ is independent of $\A_{|v|}.$ 

$\bullet$ $E(L^v_w\odot C^{v w}(t,n)\mid\A_{|v w|})=0$  

  $$E(L^v_w\odot C^{v w}(t,n)\mid \A_{|v w|})=E(A^v_w(t,n)E(C^{v w}((B^v_w,J^v_w)(t,n))\mid \A_{|v w|}))=0
  $$

$\bullet$ $E(L^v_w C^{v w}(t,n) L^v_{\overline w}C^{v\overline w}(t,n))=0$ 
for $w\ne \overline w.$  

For $|w|< |\overline w|$ (or vice versa) take the conditional expectation with respect to 
$A_{|v\overline w|}$ and use the previous statement. For $|w|=|\overline w|$ use the independence of 
$C^{v w}$ and $C^{v\overline w}$ given $\A_{|v w|}$ and argue the left hand side is 
  \begin{eqnarray*}
  &&=E(A^v_w(t,n)A^v_{\overline w}(t,n)E(C^{v w}((B^v_w,J^v_w)(t,n))
  C^{v\overline w}((B^v_{\overline w},J^v_{\overline w}(t,n))\mid\A_{|v w|})\\
  &&=E(A^v_w(t,n)A^v_{\overline w}(t,n)(E(C^{v w}((B^v_w,J^v_w)(t,n))\mid\A_{|v w|}))
  (E(C^{v \overline w}((B^v_{\overline w},J^v_{\overline w})(t,n))\mid\A_{|v w|})))\\
  &&=0
  \end{eqnarray*}

$\bullet$ $E(S_m^v S_{\overline m}^v)=0$ for $m\ne \overline m.$ 

Square out and use previous results. 
  
$\bullet$ $E(S_m^v S_m^v)=E\sum_{v\in V_m}(L^v_w\odot C^{vw})^2$ 
  
Square out and use previous results.  

$\bullet$ Let $b(n,i)=\frac{i-1}{n}\vee \frac{n-i}{n}.$
Then $E(b(n,I_n))^2\le \frac{2}{3}.$ 

Notice $EI_n=\frac{n(n+1)}{2}$ and $E(I_n)^2=\frac{n(n+1)(2n+1)}{6}.$ 
  \begin{eqnarray*}
  E(b(n,I_n))^2 &=& \frac{1}{n^2}(n^2+2E(I_n)^2 -2nEI_n-2EI_n+1)\\
  &=&\ldots= \frac{2}{3}-\frac{1}{n}+\frac{1}{3n^2}\le \frac{2}{3}
  \end{eqnarray*}
  
$\bullet$ $\sup_n\sum_{w\in V_m}E(\sup_t A^v_w(t,n))^2\le \left(\frac{2}{3}\right)^m$

Let $A^v_w(*,n)=\sup_t A^v_w(t,n).$ 
Notice $\sup_i\sup_t A_i^v(t,n)\le b(n,I_n^v).$ The recursion for $A$ is 
  $$A^v_{i w}(t,n)=A^v_i(t,n)A^{vi}_w((B^v_i,J^v_i)(t,n)).$$
We obtain 
  $$\sup_t\sum_i A^v_{i w}(t,n)\le \sup_i A^v_{i w}(*,n)\le b(n,I_n^v)\sup_i A^{vi}_w(*,J_i^v)$$
This provides 
  \begin{eqnarray*}
  E((A^v_{i w}(*,n))^2) &\le& E((b(n,I_n^v))^2\sup_i 
  E((A^{vi}_w(*,J^v_i(*,n)))^2\mid \A_{|v|}))\\
  &\le& \frac{2}{3} \sup_i E((A^{vi}_w(*,J^v_i(*,n)))^2)
  \end{eqnarray*}
By an induction on the length of $w$ we obtain the claim.   

$\bullet$ $\sup_{(t,n)}E((S_m^v(t,n))^2)\le a \left(\frac{2}{3}\right)^m$

  \begin{eqnarray*}
  \mbox{ l.h.s. }&=&\sum_{w\in V_m}E ((A^v_w(t,n) C^v_w((B^v_w,J^v_w)(t,n)))^2)\\
  &=&\sum_{w\in V_m}E((A^v_w(t,n))^2E((C^v_w((B^v_w,J^v_w)(t,n)))^2\mid \A_{|v w|}))\\
  &\le & a\sum_{w\in V_m}E((A^v_w(t,n))^2)\\
  &\le& a\sup_{t}\sum_{w\in V_m}E((A^v_w(t,n))^2)\\
  &\le& a\sum_{w\in V_m}\sup_t E((A^v_w(t,n))^2)\\
  &\le& a\left(\frac{2}{3}\right)^m
  \end{eqnarray*}

$\bullet$ $R^v_m(t,n)$ is an $L_2$-martingale in $m$ with respect to $\A_{|v|+m}$ 
for all $(t,n)$ and $v\in V.$ 

The martingale property follows by $E(S^v_m(t,n)\mid \A_{|v|+m})=0,$   
the $L_2$-statement by the previous claim. 

$\bullet$ $\sup_{t,n}E((R^v-R^v_m)(t,n))^2\le \frac{3a}{2}\left(\frac{2}{3}\right)^m$

  $$E((R^v-R^v_m)(t,n))^2\le \sum_{i\ge m}E(S^v_i(t,n))^2\le a \sum_{i\ge m}\left(\frac{2}{3}\right)^i
    =\frac{3a}{2}\left(\frac{2}{3}\right)^m$$

$\bullet$ $E(R^v_m(t,n)-R^v_m(t,\infty))^2\to_{n\to\infty}0$
for all $t\in[0,1]$ and $m\in\N.$ 

Without loss of generality let $t$ be none of the finitely many splitting points of 
the tree $v V$ up to depths $m$. Estimate $\|R^v_m(t,n)-R^v_m(t,\infty)\|_2$ by the finite sum over all 
$w\in V_{<m}$ of the terms $\|L^v_w\odot C^{v w}(t,n)-L^v_w\odot C^{v w}(t,\infty)\|_2.$
We shall show every such term converges to $0.$ 
  \begin{eqnarray*}
  \mbox{ l.h.s. }&=&  
    \|A_w^v(t,n)C^{v w}((B_w^v,J_w^v)(t,n))-A_w^n(t,\infty)C^{v w}((B_w^v,J_w^v)(t,\infty))\|_2\\
  &\le& \|(A_w^v(t,n)-A_w^v(t,\infty)) C^{v w}((B_w^v,J_w^v)(t,n)) \|_2\\
  && + \|A_w^v(t,\infty)(C^{v w}((B_w^v,J_w^v)(t,n))-C^{v w}((B_w^v,J_w^v)(t,\infty))\|_2\\
  &\le & \|A_w^v(t,n)-A_w^v(t,\infty)\|_2 \sup _{s\in[0,1]}\sup_i\|C^{v w}(s,i)\|_2\\
  && + \|C^{v w}((B_w^v,J_w^v)(t,n))-C^{v w}((B_w^v,J_w^v)(t,\infty)\|_2\\
  \end{eqnarray*}
The first term converges as $n\to\infty$ to $0$ since the difference converges to $0$ 
a.e. and is uniformly bounded by $1.$ 

Estimate the second term by the triangle inequality
  $$\le \|C(U^{v w}_{J_w^v(t,n)})-C(U^{v w})\|_2
    +\|\1_{B^v_w(t,n)\ge U^{v w}_{J^v_w(t,n)}}
    J^v_{w1}(t,n)Q^{v w 1}_{J^v_{w1}(t,n)}
    -\1_{B^v_w(t,\infty)< U^{v w}}U^{v w}Q^{v w 1}\|_2$$
The first term will converge to $0.$ Argue $U_m^{v w}$ converges a.e. to $U^{v w}$ for $m\to\infty.$ Then dominated 
convergence provides the statement, since $J_w^v(t,n)\to_n\infty$ a.e. and the function 
$C$ is bounded. 

Estimate the second term by 
  \begin{eqnarray*}
  &\le& \|(\1_{B^v_w(t,n)\ge U^{v w}_{J^v_w(t,n)}}-\1_{B^v_w(t,\infty)< U^{v w}})
    J^v_{w1}(t,n)Q^{v w 1}_{J^v_{w1}(t,n)}\|_2\\
  &&+  \|\1_{B^v_w(t,\infty)< U^{v w}}(J^v_{w1}(t,n)-U^{v w}) Q^{v w 1}_{J^v_{w1}(t,n)}\|_2\\ 
  &&+  \|\1_{B^v_w(t,\infty)< U^{v w}}U^{v w}( Q^{v w 1}_{J^v_{w1}(t,n)}-Q^{v w 1})\|_2\\
  &\le& \|\1_{B^v_w(t,n)\ge U^{v w}_{J^v_w(t,n)}}-\1_{B^v_w(t,\infty)< U^{v w}}\|_2\sup_m\|Q_m\|_2  \\
  &&+  \|J^v_{w1}(t,n)-U^{v w}\|_2\sup_m\|Q_m\|_2\\ 
  &&+  \|Q^{v w 1}_{J^v_{w1}(t,n)}-Q^{v w 1}\|_2\\
  \end{eqnarray*}
The first term converges to $0$ since the difference is bounded and converges a.e. to $0.$ The second term 
converges to $0$ since the difference is bounded and converges a.e. to $0.$ For the third term notice 
$J^v_{w1}(t,n)$ converges a.e. to $\infty$ and $b_m:=\|Q_m-Q\|_2\to_{m\to\infty}0.$ Argue 
  \begin{eqnarray*}
  \|Q^{v w 1}_{J^v_{w1}(t,n)}-Q^{v w 1}\|^2_2&=& 
    E(E((Q^{v w 1}_{J^v_{w1}(t,n)}-Q^{v w 1})^2\mid \A_{|v w 1|}))\\
  &=&\|b_{J^v_{w1}(t,n)}\|_2^2  \to_n0
  \end{eqnarray*}




Combining the above results we obtain the Theorem.\hfill q.e.d.

We come now to the weak convergence of the processes. For a vector 
$\underline t=(t_1,t_2,\ldots,t_k)\in T^*,\ k\in\N$ and a real valued function $f:T\to \R$  
let $f(\underline t)$ be the vector $(f(t_1),f(t_2),\ldots,f(t_k))$. A finite dimensional
distribution of a process $X=(X(t))_{t\in T}$ is the distribution of 
$X(\underline t),\underline t\in T^*.$ A process $X_n$ converges weakly to a process $X$,
if all finite dimensional distributions converge, i.e. $X_n(\underline t)$ converges in 
distribution to $X(\underline t)$ for all $\underline t\in T^*.$ 

  \begin{Corollary}\label{finite dimensional Convergence}
  The process $(Y_n(t))_{0<t\le 1}$ converges weakly to the process $(Y(t))_{0<1\le 1}.$  
  \end{Corollary}
This is an immediate consequence of Theorem \ref{4.8}. \hfill q.e.d. 

The convergence is stated for the half open interval $(0,1].$ We could also obtain convergence on $[0,1]$ 
by redefining $Y_n$ appropriate. However then we should be very careful about the recurrence relation. 
See the remark on left or right continuity (on $D$ and $E$) at the end of the introduction. 

  \begin{Proposition}\label{5.5}
  If $\frac{l_n}{n}\to_n t\in(0,1],\ \frac{i_n}{n}\to_n x$ and $t\ne x$ then 
    $$C(n,l_n,i_n)\to_n C(t,x)$$ 
  where the $C$-functions are given in (\ref{1.3a}) and $(\ref{1.5}).$   
  \end{Proposition}
Proof: The recursion (\ref{1.2}) for $n\ge 2$ and $1\le l\le n$ implies a recursion for 
$a(n,l)=EX(n,l)$ in $n\ge 2.$ The solution is given in (\ref{1.2a}) \cite{Mar}. 
The solution provides the correct values $EX(n,l)$ for the $(n,l)$-tuples  $(1,1),(1,0),(0,0)$
but not for the tuples $(n,0)$ for $n\ge 2.$ Therefore we have to be careful by plugging in.

The asymptotics of the harmonic numbers are  
  $$H_n=\ln n+\gamma+\frac{1}{2n}-\frac{1}{12n^2}+\frac{10}{120n^4}+O(n^{-6})$$
with $\gamma=0,577215...$ the Euler constant. We will use $H_n=\ln n+\gamma +b_n$ 
with $b_n=O(\frac{1}{n}).$ 

For $n\ge 2$ and $1\le l=l_n,i=i_n\le n$ argue  
  \begin{eqnarray*}
  C(n,l,i)&=&I+II+III\\
  I&=& \1_{l<i}\frac{1}{n}(n-1+a(i-1,l)-a(n,l))\\
  II&=&\1_{l=i}\frac{1}{n}(n-1+a(i-1,i-1)-a(n,l))\\
  III&=&\1_{l>i}\frac{1}{n}(n-1+a(i-1,i-1)+a(n-i,l-i)-a(n,l))\\
  \lim_{n\to\infty}II&=&0\\
  \lim_{n\to\infty}I&=&1_{t<x}(1+\lim_n\frac{1}{n}(2(i-1)-6(l)-6-2n+6(l)-6\\
  &&+2iH_{i-1}-2(i-l+2)H_{i-l}-2(n+1)H_n+2(n-l+3)H_{n-l+1}))\\
  &=&\1_{t<x}(-1+2x+\lim_n\frac{2}{n}(i(\ln (i-1)+\gamma+b_{i-1})\\
  &&-(i-l+2)(\ln(i-l)+\gamma+b_{i-l})\\
  &&  -(n+1)(\ln n+\gamma+b_n)+(n-l+3)(\ln(n-l+1)+\gamma+b_{n-l+1}))))\\ 
  &=&\1_{t<x}(-1+2x+\lim_n\frac{2}{n}(i\ln \frac{i-1}{n}-(i-l+2)\ln\frac{i-l}{n}\\
  && -(n+1)\ln\frac{n}{n}+(n-l+3)\ln\frac{n-l+1}{n}))\\ 
  &=&\1_{t<x}(-1+2x+2x\ln x-2(x-t)\ln (x-t)+2(1-t)\ln (1-t))\\ 
  \lim_{n\to\infty}III&=& \1_{t>x}(1+\lim_n\frac{1}{n}(2(i-1)-6(i-2)+2(n-i)-6(l-i-1)-2n+6l)\\
  &&+2(n-i+1)H_{n-i}-2(n-l+3)H_{n-l+1}+2iH_{i-1}-6H_1\\
  &&-2(n+1)H_n+2(n-l+3)H_{n-l+1}))\\
  &=&\1_{t>x}(1+\lim_n\frac{2}{n}((n-i+1)(\ln(n-i)+\gamma+b_{n-i})\\
  &&-(n-l+3)(\ln (n-l+1)+\gamma+b_{n-l+1}) +i(\ln i+\gamma+b_i)\\
  && -(n+1)(\ln n+\gamma+b_n)+ (n-l+3)
   (\ln( n-l+1)+\gamma+b_{n-l+1})))\\
  &=& \1_{t>x}(1+ 2(1-x)\ln(1-x)-2(1-t)\ln (1-t)+2x\ln x + 2(1-t)\ln (1-t)).
  \end{eqnarray*}
The statement follows easily.

\hfill q.e.d.

\end{document}